

\input amstex
\documentstyle{amsppt}
\magnification\magstep1\hsize6truein\vsize8.5truein
\rightheadtext\nofrills{Automorphisms and the K\"ahler cone of certain
Calabi-Yau manifolds}
\leftheadtext\nofrills{A. Grassi and D. R. Morrison}
\topmatter
\title
Automorphisms and the K\"ahler cone of certain Calabi-Yau manifolds
\endtitle
\author Antonella Grassi and David R. Morrison \endauthor
\address A.G.:  Mathematical Sciences
Research Institute, 1000 Centennial Drive, Berkeley, CA 94720.
 \endaddress
\address D.R.M.: Department of Mathematics, Duke University, Box 90320,
Durham, NC
27710; {\it Current address:} School of Mathematics, Institute for
Advanced Study, Princeton, NJ 08540. \endaddress
\endtopmatter
\def\NE{\operatorname{NE}}
\def\Pic{\operatorname{Pic}}
\def\Aut{\operatorname{Aut}}

\def\bkx{\overline { \Cal K (X)}}

\def\bks{\overline { \Cal K (S)}}

\def\bksp{\overline { \Cal K (S_1)}}

\def\bksq{\overline { \Cal K (S_2)}}

\def\n{\overline {   \NE(X)}}
\def\ns{\overline {   \NE(S)}}

\def\P{\Bbb P ^1}

\def\ccdot{\cdot}

\document

Let  $X$  be the fiber product over $\P$ of two rational elliptic
surfaces as in the diagram:
$$
\alignat3
& &&\ \ X &&\\
 & ^{p_1}  \swarrow && &&\searrow ^{p_2}\\
& S_1 &&  &&\  \ \ \ S_2\\
& _{\pi _1} \searrow && && \swarrow _{\pi _2}\\
& &&\ \ \P &&
\endalignat
$$
and let $f =  {\pi _1} \ccdot p =  {\pi _2} \ccdot p_2 : X \rightarrow \P$.
Schoen [S] has shown that if the surfaces are sufficiently general,
then $X$ is a smooth
Calabi-Yau manifold.
Let $\bkx$ be the closure of the K\"ahler cone of $X$.
By construction $\bkx$ has infinitely many edges. Here
 we show that there is a fundamental domain which is a (finite) rational
polyhedral cone, for the induced action of $\Aut(X)$ on $\bkx$.

Our work was inspired by some recent conjectures of the second author
[M1, M2] which derive from the ``mirror symmetry'' phenomenon for Calabi-Yau
manifolds.  In [M1], some of the data from the topological field theories
introduced by Witten [W1, W2] is used to construct some novel variations
of Hodge structure from Calabi-Yau manifolds.  In [M2], the implications
of this construction for possible compactifications of moduli spaces are
explored.  In particular, it is pointed out there that Looijenga's semi-toric
compactification method [L2] can be fruitfully applied in this situation
provided that the action of the fundamental group on the
closure of the K\"ahler cone
has a rational polyhedral fundamental domain.  This paper provides the
first non-trivial example of such a structure.

General results of Wilson [Wi] tell us that away from its intersection with
the cubic cone $W^*$ defined by cup-product, the closure $\overline{\Cal K}$
of the K\"ahler cone of a Calabi-Yau threefold is locally rational
polyhedral, with codimension one faces corresponding to primitive
birational contractions.  The number of such faces may be infinite,
and they may accumulate towards $W^*$.  The conjectures in [M2] and
the example here suggest that the automorphism group of the threefold
should fully account for such infinite phenomena in the cone.

Several authors have considered related questions for other examples.
Borcea [Bo] studied desingularized Horrocks-Mumford quintics $\widetilde{V}$
which have a K\"ahler cone whose closure has infinitely many edges.
One of the edges
is an accumulation point, and he shows that the others form a single
$\Aut(\widetilde{V})$-orbit.  Oguiso [Og] has also
constructed Calabi-Yau manifolds $Y$
for which the
closure of the K\"ahler cone has infinitely many edges, and has asked
whether the number of $\Aut(Y)$-orbits of the set of all fiber space
structures on $Y$ is finite.
 The finiteness
of the sets of $\Aut(X)$-orbits of edges and of fiber space structures
 for the examples $X$ considered here is an obvious corollary
of our main result.

 Roughly speaking the proof of our theorem goes as follows.
We first analyze
the action of the automorphism group of a general rational elliptic
surface $S$ on the cone $\bks$. Then we
describe  $\bkx$ in term of the subcones $\bksp$ and $\bksq$.  We use
several results of Looijenga [L1] about the Weyl group associated to an
elliptic ruled surface.
Another key fact  (due to Namikawa [N]) is
that  $\Aut(X)= \Aut (S_1) \times _{\P} \Aut(S_2)$.

\bigpagebreak

{\bf \S 1}

We start with some background results from
[L1].
By $\pi: S \rightarrow \P$ we denote a  general rational elliptic surface.
Any such surface is obtained  by blowing up $\Bbb P ^2$ at $9$ points.
Let $e_1, \cdots, e_9$ be the  corresponding exceptional divisors and
$h$ be the pullback of the  ample divisor on $\Bbb P^2$.
Then all fibers of $\pi$ are irreducible  and we let
$f = - K_S = 3h - \sum e_i $
denote a general fiber of the fibration.
 $S_\eta$ denotes the
generic fiber, which is an elliptic curve over the function field of
$\P$.

We choose  $e_1$ to be the ``zero-section'' of our elliptic fibration.
With that choice, any other section $e$ determines an element in
$\Pic^0(S_\eta)$ as follows.  Regarding $e$ and $e_1$ as divisors of degree
1 on $S_\eta$,  their difference is a divisor of degree 0.
We let $[e-e_1]\in\Pic^0(S_\eta)$ denote the divisor class.

Manin [Ma] has shown that addition in $\Pic^0(S_\eta)$ is related
to linear equivalence
on $S$ as follows.  Given integers
$a_2$, \dots, $a_9$, the element
$$- \sum_{i=2}^9 a_i [e_i-e_1]\in\Pic^0(S_\eta)$$
can be written in the form $[e_a-e_1]$ for some section $e_a$.
Then the linear equivalence class of $e_a$ is
  $$  3d h -(d-s-1) e_1 - (d + a_2) e_2 - \dots - (d + a_9) e_9 \in\Pic(S) $$
(where this time the addition is of divisors on $S$,
rather than on $S_\eta$), and where
   $$  d = \sum_{i=1}^9 (a_i)^2
              + \sum_{2 \le j < k \le 9}  a_j a_k
              + \sum_{i=1}^9  a_i                     $$
and
   $$  s = \sum_{i=1}^9  a_i .                         $$

The translations on  $S_\eta$ extend to automorphisms of $S$
(see for example
[MP]).
If we  choose $e_1,\dots,e_9$ so that
$[e_j -e_1]$ is of infinite order in $ S_\eta$, then the
successive translations by $e_j$
 give rise to infinitely many sections on $S$.
These sections $\{e\}$ are rational curves such that $e^2 = -1$.
The  cone of effective curves (also called the Mori cone)
 $\n$ thus has infinitely many edges. In fact,
the set of edges of $\n$ consists exactly  of the fiber $f$ and all
the sections $\{e\}$.

\proclaim{Theorem 1.1 (Manin [Ma, Theorem 6])}
Let  $\Cal T$ be the translation group  generated by all sections and  let
 $\Cal T _0$ be the
subgroup generated by the translations by $e_1, \cdots, e_9$.
Then $\Cal T_0$ is of index three in $\Cal T$ and $\Cal T$ is generated by
$\Cal T_0$ together with $1/3 \sum _{i=2} ^9 [e_i - e_1]$. (The sum
 is in the group law of $S_\eta$.)
\endproclaim

\definition{Definition-Proposition
 1.2}
 The \it root lattice \rm is the orthogonal complement of $f$ in $\Pic(S)$.
This root lattice has an integral basis
  $B = \{ e_1 - e_2,  e_2 -e_3, \cdots ,
e_8 - e_9, h-e_1- e_2 -e_3 \}$,
determined by a realization of $S$ as the blowup of ${\Bbb P}^2$ in 9 points.
\enddefinition

\definition{Definition 1.3}  For any element $\alpha \in B$, define the
\it fundamental reflection  $s_\alpha$ \rm
to be the automorphism of $\Pic(S)$ given by:
$$s_\alpha(x) = x + (x \cdot \alpha) \alpha .$$
 The group $W$ generated by these reflections is called
the \it Weyl group \rm of $B$.
\enddefinition

\definition{Definition 1.4}

(1.4.1) $\Cal C = \{ x \in \Pic (S) \text{ such that }
x \cdot \alpha > 0 , \ \forall \alpha \in B \}$
is called the
\it fundamental chamber.\rm

(1.4.2)  The $W$-orbit of $\overline {\Cal C}$ is called the
\it Tits cone.\rm
\enddefinition

\proclaim{Theorem 1.5 (Bourbaki, [B, Ch.~V, \S 4])}
  The fundamental chamber
   $\overline {\Cal C}$ is a fundamental domain for the induced
action of $W$ on the Tits cone, that is, the $W$-translates of
$\overline {\Cal C}$ define a partition of the Tits cone.
\endproclaim

\vskip 0.2in

{\bf \S 2 The surface case.}

\bigpagebreak

The following two theorems relate $W$ and the Tits cone with $\bks$
and the sections.

\proclaim{Theorem 2.1 (Looijenga [L1, \S 3 - \S 4])}

(2.1.1)  Any $\{w(e_j) \}$,
for $w \in W$ is an exceptional curve .
Conversely, if $\{e'_j \}$ is an exceptional curve,
then  $e'_j = w(e_j)$ for some $w \in W$.

(2.1.2) $\bks $ is contained in the Tits cone.
\endproclaim

\proclaim{Theorem 2.2} $\Cal T \subset W$.
\endproclaim

\demo{Proof} By Manin's Theorem 1.1 it is enough to check the statement
for the translation by  each $e_i$ and
the translation by  $1/3 \oplus _{i=2} ^9 e_i$.
The idea of the proof is the following:
take any $x \in {\Cal C}$ and consider $t x$, where $t \in \Cal T$.
 By 1.5 we know that there exists an element $w \in W$ such that
$w tx \in {\Cal C}$.
One wants to show that $wt$ is the identity map.

This procedure can be implemented explicitly by
following Bourbaki's argument:  successively
reflect $tx$ with respect to the
walls (of $\bar{\Cal C}$) between $tx$ and ${\Cal C}$,
producing $w\in W$ such that
$wtx\in {\Cal C}$.  In each case, we carried this out,
and then verified that $wt$ is the identity map using {\smc maple}.

We illustrate the form of the results of our computation with an example.
Suppose that we wish to write the translation $t_2$ with respect to
the element $e_2$ as a product
 of reflections. Following Manin we choose $h, -e_1, \cdots,- e_9$ as a
basis of $\Pic(S)$ and write
any element in $\Pic(S)$ as $x= [b, a_1, \cdots, a_9]$. Note that
every permutation
 of the elements $a_1, \cdots ,a_9$
can be realized as a composition of reflections; we use cycle notation
for such permutations.
Then our calculations give:
$$P_6 w_s P_5 w_s P_4 w_s P_3 w_s P_2 w_s P_1 t_2 = Id, \text{ where }$$
$$ \aligned
P_6 &= (2,3,4,1,5,6,7,8,9) \\
P_5 &=(1,5,6,2,3,4,7,8,9) \\
P_4 &=(7,8,9,1,2,3,4,5,6) \\
P_3 &= (4,5,6,1,2,3,7,8,9) \\
P_2 &= (1,7,8,2,3,4,5,6,9) \\
 P_1 &= (1,9,2,3,4,5,6,7,8)
\endaligned$$
and $w_s$ denotes the reflection with respect to the root $h-e_1 - e_2- e_3$.
\qed
\enddemo

\proclaim{Theorem 2.3} Let $S$ be a general rational elliptic surface.

 (2.3.1)  There exists a fundamental domain  $\overline {\Cal D}$ for
the induced action of the automorphism group  of $S$ on the Tits cone.

(2.3.2) Furthermore   $\bks \cap \overline {\Cal C}$ and
$\bks \cap \overline {\Cal D}$ are (finite) rational polyhedral cones.
\endproclaim

\demo{Proof}

 (2.3.1) It  is  sufficient to show that  there exists a fundamental
domain for the induced action of the translation group $\Cal T$ on the
Tits cone. In fact,
$\Cal T$ is a subgroup of the automorphism group of $S$
of finite index.  This follows from the facts that any automorphism
of $S$ can be translated to one which preserves the zero-section,
and that the group of automorphisms which preserve the zero-section
can be identified with $\Aut(S_\eta)$,
which is finite.

Let $x$  be in the Tits cone. Then  there exists an element $w \in W$
such that
$w ^{-1}(x) \in \overline {\Cal C}$ (see 1.5).
By 2.1  $w(e_9)$ is an exceptional curve and  we can find
 $T \in \Cal T$ such that $Tw(e_9) = e_9$.

$Tw$ is an orthogonal transformation, so
$Tw (e_9 ^{\perp} )\subset   e_9 ^{\perp}$.
Moreover, $Tw(f^\perp)\subset f^\perp$, so $Tw$ preserves the space
$(e_9,f)^\perp$.  Now $(e_9,f)^\perp$ is the space spanned by the
sub-root system $E_8\subset E_9$.  (This is the root system
obtained by removing the simple root $\{e_8 - e_9 \}$ from $B$.)
$Tw_{(e_9,f)^\perp}$ acts as a reflection on this space $(e_9,f)^\perp$
and, since $Tw$ lies in $W$, the restriction $Tw_{(e_9,f)^\perp}$
must actually belong to the Weyl group $W(E_8)$
(because any reflection is generated
by reflections in simple roots (Humphreys [H, \S 1.5])).

Thus, there exists an element $w' \in W(E_8)$ such that:
$$
\left\{
\aligned w' \ccdot T \ccdot w (e_9)  &= e_9    \\
w' \ccdot T \ccdot w (h -e_1 -e_2 - e_3 )&=  h -e_1 -e_2 - e_3 \\
 w' \ccdot T \ccdot w (e_i - e_{i+1})  &= e_i - e_{i+1} \quad \text{ if }
 i \neq 8 \endaligned
\right.
$$
  Hence:
 $$
\left\{
\aligned w' \ccdot T \ccdot w (e_9)  &= e_9    \\
  w' \ccdot T \ccdot w (e_i)  &= e_i  + \beta g \quad \text{ if } i \neq 9 \\
 w' \ccdot T \ccdot w (h)  &= h + 3 \beta g \endaligned
\right.
$$
for some divisor $g$ and some $\beta \in \Bbb Z$.

Note that $w'\ccdot T \ccdot w $ fixes $e_9$ and $f$.
Intersecting  both sides of the above equation  with $f$  and $e_i$
we obtain  $\beta =0$.
Thus $ w' \ccdot T \ccdot w$ is the identity map, in particular
$T^{-1} \ccdot {{w'} ^{-1}} (x)$ belongs to the fundamental chamber $\Cal C$.

This shows that $\Cal D $, the convex hull of $ \cup _{w'_i \in W(E_8)}
\overline {w'_i(\Cal C)}$ is a fundamental domain
for the induced action of $\Cal T$ on $\Pic(S)$.

(2.3.2) By  2.1 $\bks$ is contained in the Tits cone of $S$.
 It is enough to show that $\bks$ cuts the fundamental chamber
$\overline {\Cal C}$
in a finite rational
polyhedral cone.  In fact $W(E_8)$ is a finite group and thus
$\overline {\Cal D}$ is the union of finitely many $W$-translates of
$\overline {\Cal C}$.

Note that the walls of $\overline {\Cal C}$ in $\Pic(S)$  are defined
by the hyperplanes
corresponding to the roots: $\{ e_1 - e_2,  e_2 -e_3, \cdots ,
e_8 - e_9, h-e_1- e_2 -e_3 \}$.
The statement follows by duality from the following Lemma 2.4. \qed
\enddemo

\proclaim{Lemma 2.4} Let $\ns$ be the cone of effective curves on $S$.
Then the convex hull of $\ns$ and $B=\{ e_1 - e_2,  e_2 -e_3, \cdots ,
 e_8 - e_9, h-e_1- e_2 -e_3 \}$ is a finite rational polyhedral cone.
\endproclaim

\demo{Proof}
We claim that the convex hull of $\ns$ and $B$ is
the union of finitely many translates of the finite rational
polyhedral cone spanned by $B$ and $\{e_1,\dots,e_9\}$.

The fiber $f=3h-\sum e_i$ is easily seen to be a convex combination of
elements of $B$.  Every other edge of $\ns$ corresponds to a section
$\sigma$ of the elliptic fibration, and each such $\sigma$ has the
form $\sigma=T(e_1)$ for some $T\in\Cal T$.
Since $\Cal T_0\subset\Cal T$ is of finite index (Manin's theorem 1.1),
it will suffice to show that each $\sigma=T(e_1)$ with $T\in\Cal T_0$
is in the finite rational polyhedron spanned by $B$ and $\{e_1,\dots,e_9\}$
(i.e., that it is a convex combination of $B$ and $\{e_1,\dots,e_9\}$).

{}From \S 1 , we have:

$$
\sigma  -e_1 =   3d h + (d-s) e_1 + (d + a_2) e_2 + \dots + (d + a_9) e_9 $$
Thus
$$ \aligned
\sigma  -e_1  &= 3d (h-e_1 -e_2 -e_3)\\
 &+ (2d  +s)(e_1-e_2)\\
&+
(4d+  s - a_2)(e_2 -e_3) \\
&+ (6d + s - a _2 - a_3)( e_3 -e_4) \\
&+ (5 d + s - a _2 - a_3- a_4)(e_4 -e_5) \\
&+ (4d  + a _6 + a_7+ a_8 +a_9) (e_5 -e_6)\\
&+ (3d  + a_7+ a_8 +a_9) (e_6 -e_7)\\
&+  (2d +   a_8 +a_9) (e_7 -e_8)\\
&+   (d +   a_9) (e_8 -e_9)
\endaligned$$

\smallpagebreak

Our goal is to show that all the above coefficients are non negative.
If
$\sigma \neq  e_2, \cdots e_9$, then $(\sigma-e_1) \cdot e_j \geq 0$,
for $j \neq 1$,
and also  $(\sigma - e_1) \cdot (h - e_2 - e_3) \geq 0$.
   Thus all the coefficients  in the above formula are non negative
except $2d  +s$.
By Manin's formula we have
$$\aligned
2d +s &= 2 \sum (a_i)^2 + 2 \sum _{j<k} a_j a_k + 3 \sum a_i\\
 &= \sum (a_i)^2  + ({\sum a_i + \frac{3}{2}})^2 -\frac{9}{4}
\endaligned$$

Thus $2d +s = \sum (a_i)^2  + ({\sum a_i + \frac{3}{2}})^2 -\frac{9}{4} <0$
only if $a_i=-1$, for some index $i$ and $a_j = 0$, for $j \neq i$.
But then $\sigma=e_i$.
\qed
\enddemo

\vskip 0.2in

 {\bf \S 3 The threefold case.}

\bigpagebreak

 Let $ F$ be a general divisor in the fiber of
$ f : X \rightarrow \P$. Denote by $f_i$  a fiber of the
fibration $S_i \rightarrow \P$.
Then  $ F = p _i ^*(f_i)$.

\proclaim{Proposition 3.1} $\bkx = \{  p_1^*( L) \otimes p_2^*( M),
\text{ such that }
\  L \in \bksp \text{ and }  M \in \bksq \}$.
(The decomposition of elements as $p_1^*( L) \otimes p_2^*( M)$
is not necessarily unique.)
\endproclaim

\demo{Proof}

Let $A$ be a nef and big divisor on $X$. Then we can write $A =  p^*( A_1)
\otimes p_2^*( A_2)$, for some divisors
 $A_i \in \Pic(S_i)$; following Namikawa [N, Proof of 2.1], we can
ensure that $A_1$ and $A_2$ are effective.

Suppose that $A_1$ is not nef.
If $C$ is an irreducible curve on $S_1$ with $A_1\cdot C<0$, then $C$
must be a component of the effective divisor $A_1$
of negative self-intersection,
hence an exceptional curve of the first kind.  There will be finitely
many such curves $C_1$,\dots,$C_k$, and they are in addition sections
of the fibration $\pi_1$.  If we let $m=\max\{-A_1\cdot C_i\}$, then
$A_1+mf_1$ is a nef divisor on $S_1$ which is not ample (since its
intersection number with some $C_i$ is 0).

On the other hand, if $A_1$ \it is \rm nef, and we let
$$m=-\min\{A_1\cdot C \,|\, C \text{ is an exceptional
curve of the first kind}\},$$
then $A_1+mf_1$ is still a nef divisor, since $(A_1+mf_1)\cdot C\ge0$
for every exceptional curve $C$, and $(A_1+mf_1)\cdot f_1=A_1\cdot f_1\ge0$.
(We are using here the fact that the edges of the Mori cone are exactly
the $C$'s, and $f_1$.)  Moreover, for some exceptional curve $C$,
$A_1+mf_1\cdot C=0$.

Thus, in either case, replacing $A_1$ and $A_2$ by $A_1+mf_1$ and
$A_2-mf_2$, respectively, we may assume that $A_1$ is a nef divisor
which is not ample, i.e., that for some exceptional curve $C$ on $S_1$,
$A_1\cdot C=0$.  We claim that $A_2$ is nef.  For if $\Gamma$ is
any curve on $S_2$, consider the curve
$C\times_{\P}\Gamma$ on $X$.  Since $A$ is nef on $X$, we have
$$0\le A\cdot(C\times_{\P}\Gamma)
=A_1\cdot C+A_2\cdot\Gamma=A_2\cdot\Gamma.$$
 \qed
\enddemo

\proclaim{Corollary 3.2} $\Pic(X)$ has rank $19$.
\endproclaim

\proclaim{Proposition 3.3 (Namikawa [N])}
$\Aut(X)= \Aut (S_1) \times _{\P} \Aut(S_2)$.
\endproclaim

Combining 3.1 and 3.3  with 2.3 we have:

\proclaim{Theorem 3.4}
 There is a fundamental domain which is a (finite) rational
polyhedral cone, for the
induced action of $\Aut(X)$ on $\bkx$.
\endproclaim

\demo{Proof} Use $(\bksp\cap\overline{\Cal D_1})\times_{\Bbb P^1}
(\bksq\cap\overline{\Cal D_2})$.
\qed
\enddemo

\enddocument